\documentclass[leqno,12pt]{amsart}

\usepackage{latexsym,esint}

\theoremstyle{plain}

\def\endproof{\hspace*{\fill}\mbox{\ \rule{.1in}{.1in}}\medskip }

\newtheorem{theorem}{Theorem}[section]

\newtheorem{lemma}[theorem]{Lemma}

\theoremstyle{definition}

\numberwithin{equation}{section}

\begin{document}
\title[traveling waves in the 3D Boussinesq system]
{On the existence of traveling waves in the 3D Boussinesq system}
\author{Marta Lewicka}
\address[M. Lewicka]{University of Minnesota, Department of Mathematics, 
127 Vincent Hall, 206 Church St. S.E., Minneapolis, MN 55455, USA }
\email{lewicka@math.umn.edu}
\author{Piotr B. Mucha}
\address[P.B. Mucha]{Institute of Applied Mathematics and Mechanics, 
 University of Warsaw, ul. Banacha 2, 02097 Warszawa, Poland}
\email{p.mucha@mimuw.edu.pl}

\begin{abstract}
We extend earlier work on traveling waves in premixed flames 
in a gravitationally stratified medium, subject to the Boussinesq
approximation. For three-dimensional channels not aligned with the gravity
direction and under the Dirichlet boundary conditions in the fluid velocity,
it is shown that a non-planar traveling wave, corresponding to a non-zero
reaction, exists, under an explicit condition relating the geometry of the
crossection of the channel to the magnitude of the Prandtl and Rayleigh numbers,
or when the advection term in the flow equations is neglected.
\end{abstract}

\maketitle

\section{Introduction}

The Boussinesq-type system of reactive flows is a physical model in the
description of flame propagation in a gravitationally stratified medium
\cite{Z}.
It is given as the reaction-advection-diffusion equation for the reaction
progress $T$ (which can be interpreted as temperature), 
coupled to the fluid motion through the advection velocity, 
and the Navier-Stokes equations for the 
incompressible flow $u$ driven by the temperature-dependent force term. 
After passing to non-dimensional
variables \cite{BCR, VR1}, the Boussinesq system for flames takes the form:
\begin{eqnarray}
T_t  + u\cdot \nabla T - \Delta T &=& f(T)\nonumber\\
u_t + u\cdot \nabla u - \nu\Delta u   + \nabla p &=& T\vec\rho\label{nsb}\\
\mbox{div } u &=& 0.\nonumber
\end{eqnarray}
Here, $\nu>0$ is the Prandtl number, that is the ratio of the kinematic and thermal
diffusivities (inverse proportional to the Reynolds number).
The vector $\vec\rho = \rho\vec g$ 
corresponds to the non-dimensional gravity $\vec g$ scaled by the  Rayleigh number 
$\rho>0$. The reaction rate is given by a nonnegative 'ignition type' 
Lipschitz function $f$ of the temperature, this last one normalized to satisfy: 
$0\leq T\leq 1$.
The above model can be derived from a more complete system under the
assumption that the Lewis number equals $1$.

We study the system (\ref{nsb}) in an infinite cylinder 
$D \subset \mathbb{R}^3$ with a smooth, connected 
crossection $\Omega\subset \mathbb{R}^2$. 
Recent numerical results, motivated by the astrophysical context \cite{VR1, VR2}, 
suggest that the initial perturbation in $T$ either quenches or develops
a curved front, which eventually stabilizes and propagates as a traveling wave.
On the other hand, existence of non-planar traveling waves for
the single reaction-advection-diffusion equation in a prescribed
flow has been a subject of active study in the last decade \cite{Xin, B, R}.
For system (\ref{nsb}), existence of traveling waves 
has been considered  under the no-stress or Dirichlet 
boundary conditions in $u$, in channels of various inclinations and
dimensions \cite{CKR, BCR, TV, CLR, Le}. 

The main difference presents itself at the orientation of $D$ 
with respect to $\vec g$; when they are aligned there are no non-planar 
fronts at small Rayleigh numbers \cite{CKR}, while in the other case 
a traveling front, necessarily non-planar, is expected to exist at any range 
of parameters.  This has been rigorously proven: in \cite{BCR} for $n=2$ 
dimensional channels $D$ and under no-stress boundary conditions, 
in \cite{CLR} for $n=2$ and the more physical no-slip conditions
and in \cite{Le} for the same boundary conditions and arbitrary dimension $n$,
but for a simplified system (corresponding to the infinite Prandtl number
$\nu=\infty$) when the Navier-Stokes part of (\ref{nsb}) is
replaced by the Stokes system.

\bigskip

The purpose of this paper is to remove this last assumption, for three
dimensional channels. Namely, we will investigate the model 
supplied by the Navier-Stokes system. 
We assume that $\vec \rho$
is not parallel to the unbounded direction of $D$, which after an 
an elementary change of variables \cite{BCR} amounts to studying:
\begin{equation*}\label{dom}
D=(-\infty, \infty)\times \Omega = \left\{(x,\tilde x); 
~~ x\in\mathbb{R}, ~ \tilde x\in\Omega\right\}
\end{equation*}
and 
\begin{equation*}\label{grav}
\vec\rho\cdot e_3\neq 0.
\end{equation*}
We will prove the existence of a traveling wave solution to (\ref{nsb}):
$T(x-ct,\tilde x)$, $u(x-ct,\tilde x)$, with the speed $c$ to be determined
and under the boundary conditions:
\begin{equation}\label{boundary}
\frac{\partial T}{\partial \vec n} = 0
~~~~\mbox{ and } ~~~~ u = 0 \quad \mathrm{ on } ~~ \partial D,
\end{equation}
where $\vec n$ is the unit normal to $\partial D$.
Such a front satisfies:
\begin{eqnarray}
-c T_x - \Delta T + u\cdot \nabla T &=& f(T)\nonumber\\
-c u_x + d u\cdot\nabla u -\nu  \Delta u  + \nabla p &=& T\vec\rho\label{trav_wav}\\
\mathrm{div }~ u &=& 0.\nonumber
\end{eqnarray}
We set constant $d$ to be $0$ or $1$.  For the simplified system, when the
advection in $u$ has been neglected and $d=0$, the theorem below states
existence of a non-planar traveling wave, for any crossection $\Omega$,
Prandtl number $\nu$ and Rayleigh number $\rho$.  

For the full system when $d=1$, we need to assume the following relative 
thinness condition, involving $\nu$, $\vec\rho$, the area $|\Omega|$,
and the Poincar\'e and the Poincar\'e-Wirtinger 
constants $C_P$, $C_{PW}$ of $\Omega$:
\begin{equation}\label{small} 
\sqrt{14} \frac{C_P}{\nu\sqrt{\pi\nu}} |\Omega|^{1/2}
\left(|\vec\rho|C_{PW} + 
\left(\fint_\Omega |\vec\rho\cdot(0,\tilde x)|^2\right)^{1/2}\right) < 1.
\end{equation}
This condition is essential in our analysis and it is not clear if the 
below existence result holds without it.
Recall that $C_P$ is determined by the thinness of $\Omega$,
and hence (\ref{small}) admits domains with large area which are sufficiently thin.
Respectively, $C_{PW}$ depends on the maximum of (inner) distances 
between points in $\Omega$.
On the other hand, the quantities relating to smoothness of $\partial \Omega$
have no direct influence on (\ref{small}).

The nonlinear Lipschitz continuous function $f$ is assumed to be of ignition type:
\begin{equation*}\label{igni}
f(T) = 0 \mbox{ on } (-\infty,\theta_0]\cup [1, \infty),\qquad
f(T) > 0 \mbox{ on } (\theta_0, 1)
\end{equation*}
for some ignition temperature $\theta_0\in (0,1)$.

The following is our main result:
\begin{theorem}\label{th2}
Assume that either $d=0$ or $d=1$ and (\ref{small}) holds. 
Then there exist $c>0$, $T\in \mathcal{C}^{2,\alpha}(D)$ with $\nabla T\in L^2(D)$,
$u\in H^3\cap \mathcal{C}^{2,\alpha}(D)$,
$p\in\mathcal{C}^{1,\alpha}_{loc}(D)$ satisfying (\ref{trav_wav}) and (\ref{boundary})
together with:
\begin{equation}\label{tre.due}
\lim_{x\to\pm\infty} ||u(x,\cdot)||_{\mathcal{C}^2(\Omega)} =
\lim_{x\to\pm\infty} ||\nabla T(x,\cdot)||_{L^\infty(\Omega)} = 0.
\end{equation}
Moreover $T(D)\subset [0,1]$, $~~~\max_{x\geq 0, y\in\Omega }T(x,y)= \theta_0$, 
and there is a nonzero reaction: 
\begin{equation*}
\int_D f(T) \in (0,\infty).
\end{equation*}
The limits of $T$ satisfy:
$$\lim_{x\to +\infty} ||T(x,\cdot)||_{L^\infty(\Omega)} = 0,\qquad
\lim_{x\to -\infty} ||T(x,\cdot)-\theta_-||_{L^\infty(\Omega)} = 0$$ 
for some: $\theta_-\in(0,\theta_0]\cup \{1\}$. 
\end{theorem} 
The following sections are devoted to the proof of Theorem \ref{th2}.
In section 2 we formulate some auxiliary results, of an independent interest. 
In particular, we prove a weak version of Xie's
conjecture \cite{XieStokes} for the Stokes operator
(established in \cite{Xie} for the Laplacian). 
Based on results in \cite{LLP}, we then derive an a priori estimate valid in 
any channel $D$, whose cross-section $\Omega$ fulfills 
the geometrical constraint (\ref{small}). This allows us 
to obtain uniform bounds on the quantities involved in the fixed point
argument (Theorem \ref{cor1}) in sections 3 and 4; in particular the bounds
are independent of length of the compactified domains $R_a= [-a,a]\times\Omega$.
The set-up for the Leray-Schauder degree 
is different than in \cite{CLR, Le}: we solve the flow equations in the 
full unbounded channel $D$, while the reaction equation is solved in $R_a$.
Once the uniform bounds are established, we refer to \cite{BCR, Le} for
further details of the proofs.
In section 5 we improve a sufficient condition from \cite{Le} for the left
limit $\theta_-$ of the temperature profile $T$ obtained in Theorem \ref{th2}
to be equal to $1$.

We remark that the a priori estimates we derive do not preclude the
solutions $(T, u)$ to have arbitrary large norms. Indeed, the main chain
of estimates eventually leads to inequality (\ref{inq-main}), whose 
right hand side has a linear growth in terms of the left hand side, 
and thanks to condition (\ref{small}),
the main bound on $\|u\|_{L^\infty}$ does not restrict 
the magnitude of this quantity. 
Similar estimates are known also for 
solutions to the Navier-Stokes equations  for the $2$d and cylindrical
symmetric systems \cite{La,Yu,Mu1,Mu2}, and in presence of 
a special geometrical constraint on the domain \cite{Mu1,Mu2}.

We will always calculate all numerical constants at the leading 
order terms explicitly.
By convention, the norms of a vector field $u$ on $D$ are given as:
$\|u\|_{L^\infty(D)} =
\left(\sum_{i=1}^3\|u^i\|^2_{L^\infty(D)}\right)^{1/2}$
and $\|u\|_{L^2(D)} = \left(\sum_{i=1}^3\|u^i\|^2_{L^2(D)}\right)^{1/2}$. 

\bigskip

\noindent{\bf Acknowledgments.}
M.L. was partially supported by the NSF grant DMS-0707275
and by the Center for Nonlinear Analysis (CNA) under
the NSF grants 0405343 and 0635983. P.B.M.  
has been supported by
MNiSW grant No. N N201 268935 and by ECFP6 M. Curie ToK program SPADE2,
MTKD-CT-2004-014508 and SPB-M.

\section{Auxiliary results}

In the sequel, we will need a uniform estimate for the supremum of 
the solution to Stokes system in $D$.
The known proofs of the inequality $\|u\|_{L^\infty}\leq 
C_\Omega \|\nabla u\|_{L^2}^{1/2}\|\mathcal{P}\Delta u\|_{L^2}^{1/2}$,
$\mathcal{P}$ being the Helmholtz projection,  
are based on the a-priori
estimates in \cite{ADN}, which hold for smooth domains. Therefore the constant 
$C_\Omega$ depends strongly on the boundary curvature, and becomes unbounded as 
$\Omega$ tends to any domain with a reentrant corner.  It has been conjectured
by Xie \cite{XieStokes} that actually $C_\Omega=1/\sqrt{3\pi}$.  To our
knowledge, this is still an open question.  Below we prove its weaker version,
sufficient to our purpose and involving lower order terms.

\begin{theorem}\label{thour_bound}
Let $g\in L^2(D)$.  Then the solution $u\in H^2\cap H_0^1(D)$ 
to the Stokes system:
\begin{equation}\label{stokes}
- \nu\Delta u + \nabla p = g, \qquad
\mathrm{div }~ u = 0 \quad \mbox{ in } D
\end{equation}
satisfies the bound:
$$\|u\|_{L^\infty(D)} \leq \frac{2}{\sqrt{2\pi\nu}} 
\|\nabla u\|^{1/2}_{L^2(D)}\|g\|^{1/2}_{L^2(D)} 
+ C_\Omega\|\nabla u\|_{L^2(D)},$$
where constant $C_\Omega$ depends only on the crossection $\Omega$.
\end{theorem}
\begin{proof}
{\em 1.} We first quote two results, whose combination will yield the proof.
The first one is Xie's inequality \cite{Xie} 
for the Laplace operator in a 3d domain. 
Namely, for any $u\in H^2\cap H_0^1(D)$ there holds:
\begin{equation}\label{thXie}
\|u\|_{L^\infty(D)} \leq 
\frac{1}{\sqrt{2\pi}} \|\Delta u\|_{L^2(D)}^{1/2} \|\nabla u\|_{L^2(D)}^{1/2}.
\end{equation}
The crucial information in the above estimate is that the constant 
$1/\sqrt{2\pi}$ is good for all open subsets of $\mathbb{R}^3$.

\medskip

The next result is a recent commutator estimate by Liu, Liu and Pego \cite{LLP}. 
Recall first \cite{S} that for any vector field $u\in L^2(D)$ there exists 
the unique decomposition $u=\mathcal{P}u + \nabla q$ with 
$\mbox{div}(\mathcal{P}u) =0$, and $q$ solving in the sense of distributions:
$$\Delta q = \mbox{div }u \quad \mbox{ in } D, \qquad 
\frac{\partial q}{\partial\vec n} =0 \quad \mbox{ on } \partial D.$$
This Helmholtz projection satisfies: 
$\|\mathcal{P}u\|_{L^2(D)}\leq \|u\|_{L^2(D)}.$
In this setting, it has been proved in \cite{LLP}
that for every $\epsilon > 0$ there exists 
$C_{\epsilon,  \Omega}>0$ such that:
\begin{equation}\label{thLLP}
\forall u\in H^2\cap H_0^1(D) \qquad 
\int_D |(\Delta\mathcal{P} - \mathcal{P}\Delta) u|^2 \leq 
\left(\frac{1}{2} + \epsilon\right) \int_D |\Delta u|^2 
+ C_{\epsilon,  \Omega}\int_D |\nabla u|^2.
\end{equation}
\noindent The proof in \cite{LLP}, written for bounded domains,  
can be directly used also for the case of cylindrical domains $D$
with smooth boundary (since the covering number for the partition of unity 
on $\partial D$ is finite).

\medskip

{\em 2.}
Applying the Helmholtz decomposition to (\ref{stokes})  we arrive at:
$-\nu\mathcal{P}\Delta u = \mathcal{P}g$, which can be restated as:
$$-\Delta u = (\mathcal{P}\Delta - \Delta\mathcal{P}) u
+\frac{1}{\nu}\mathcal{P}g,$$
since $\mathcal{P}u = u$. Using (\ref{thLLP}) we obtain:
$$\|\Delta u\|_{L^2(D)}\leq \frac{3}{4} \|\Delta u\|_{L^2(D)} + C_\Omega
\|\nabla u\|_{L^2(D)} + \frac{1}{\nu}\|\mathcal{P}g\|_{L^2(D)},$$
which yields:
\begin{equation}\label{lapl}
\|\Delta u\|_{L^2(D)}\leq  \frac{4}{\nu}\|g\|_{L^2(D)}
+  C_\Omega\|\nabla u\|_{L^2(D)}.
\end{equation}
Now combining (\ref{lapl}) and (\ref{thXie}) proves the result.
\end{proof} 

\medskip

We will also need an extension result for divergence free vector fields.
Define a compactified domain $R_a = [-a,a]\times\Omega$.
\begin{theorem}\label{vtilde}
For any $a>0$ and any $\epsilon>0$ there exists a linear continuous 
extension operator $E:\mathcal{C}^{1,\alpha}(R_a) \longrightarrow 
\mathcal{C}^{1,\alpha}(D)$, such that for every
$u\in\mathcal{C}^{1,\alpha}(R_a)$
there holds:
\begin{itemize}
\item[(i)] $(Eu)_{\mid R_a} = u$,
\item[(ii)] if $\mathrm{div }~ u =0$ in $R_a$, then $\mathrm{div }~ (Eu) = 0$
in $D$,
\item[(iii)] $\|Eu\|_{L^\infty(D)} \leq (1+\epsilon) \|u\|_{L^\infty(R_a)}$.
\end{itemize}
\end{theorem}
\begin{proof}
Given a vector field $u\in\mathcal{C}^{1,\alpha}([-a,0]\times\Omega)$ we shall construct
its extension $\tilde u\in\mathcal{C}^{1,\alpha}([-a,\infty)\times\Omega)$
such that (ii) holds together with:
\begin{equation}\label{marta}
\|\tilde u\|_{L^\infty([-a,\infty)\times\Omega)} \leq (1+\epsilon) 
\|u\|_{L^\infty([-a,0)\times\Omega)}.
\end{equation}
This construction, being linear and continuous with respect to the
$\mathcal{C}^{1,\alpha}$ norm, will be enough to establish the lemma.

Fix a large $n>0$. For $x\in[-a,a/2n^2]$ and $\tilde x\in\Omega$, define the
vector $v(x,\tilde x)$ with components:
\begin{equation*}
\begin{split}
v^1(x,\tilde x) & = \left\{\begin{array}{ll}u^1(x,\tilde x) & \mbox{ for }
    x\in[-a,0]\\
\lambda_1 u^1(0,\tilde x) + \lambda_2 u^1(-nx, \tilde x) + \lambda_3
u^1(-n^2x,\tilde x)  & \mbox{ for } x\in[0,a/2n^2],\end{array}\right.\\
\mbox{for } i=2,3:& \\
v^i(x,\tilde x) & = \left\{\begin{array}{ll}u^i(x,\tilde x) & \mbox{ for }
    x\in[-a,0]\\
-n\lambda_2 u^i(-nx, \tilde x) - n^2\lambda_3
u^i(-n^2x,\tilde x)  & \mbox{ for } x\in[0,a/2n^2],\end{array}\right.
\end{split}
\end{equation*}
where:
$$\lambda_1=\frac{(1+n)(1+n^2)}{n^3}, \quad 
\lambda_2=-\frac{1+n^2}{n^2(n-1)},\quad
\lambda_3=\frac{1+n}{n^3(n-1)}.$$
Since we have: $\sum_{i=1}^3 \lambda_i =1$, $-n\lambda_2 - n^2\lambda_3 = 1$ and
$n^2 \lambda_2 + n^4\lambda_3 =1$, it follows that
$v\in\mathcal{C}^{1,\alpha}([-a,a/2n^2]\times\Omega)$. Also, by an explicit
calculation, we see that $\mbox{div }u=0$ implies $\mbox{div }v=0$.

Let now $\phi\in\mathcal{C}^\infty(\mathbb{R},[0,1])$ be a non-increasing cut-off
function such that $\phi(x)=1$ for $x<0$ and $\phi(x)=0$ for $x>a/3n^2$.
Define:
$$\tilde u(x,\tilde x) = \phi(x) v(x,\tilde x) +
\int_0^x \phi'(s) v^1(s,\tilde x)~\mbox{d}s \cdot e_1.$$
Clearly, $\tilde u\in\mathcal{C}^{1,\alpha}([-a,\infty)\times\Omega)$
and $\mbox{div }\tilde u=0$ if $\mbox{div } u=0$.  
Further:
\begin{equation*}
\begin{split}
|\tilde u^1(x,\tilde x)| & \leq (|\lambda_2| + |\lambda_3|)\|u^1\|_{L^\infty}
+ |\lambda_1|\cdot\left|\phi(x) u^1(0,\tilde x) 
+(1-\phi(x)) \|v^1\|_{L^\infty}\right| \\
& \leq (|\lambda_2| + |\lambda_3|)\|u^1\|_{L^\infty} +
|\lambda_1|(|\lambda_1| + |\lambda_2| + |\lambda_3|)\|u^1\|_{L^\infty},\\
|\tilde u^i(x,\tilde x)| & \leq (n|\lambda_2| + n^2|\lambda_3|)\|u^i\|_{L^\infty}
\qquad \mbox{for } i=2,3.
\end{split}
\end{equation*}
Since $\lambda_1\to 1$,  $|\lambda_2|,\lambda_3\to 0$,
$n|\lambda_2|\to 1$ and   $n^2\lambda_3\to 0$ as $n\to\infty$,
the estimate (\ref{marta}) holds if only $n$ is sufficiently large,
which ends the proof.

We remark that the norm of the operator $E$ blows up when $\epsilon\to 0$. 
Indeed, one cannot have $\epsilon = 0$ in (iii) and keep the norm
of $E$ bounded.
\end{proof}

The following elementary fact will be often used:
\begin{lemma}\label{lem_pomoc}
For any $u\in L^2\cap\mathcal{C}^{0,\alpha}(D)$ there holds:
$\lim_{x\to\pm \infty} \|u(x,\cdot)\|_{L^\infty(\Omega)} = 0.$
\end{lemma}

\section{The bound on $\|u\|_{L^\infty(D)}$}

In this section, given $c\in\mathbb{R}$, $\tau\in [0,1]$, 
a divergence-free vector field $\tilde v\in\mathcal{C}^{1,\alpha}(D)$
and a boundedly supported $\tilde T\in\mathcal{C}^{1,\alpha}(D)$, 
we consider the following problem:
\begin{eqnarray}
-cu_x - \nu\Delta u + \tau d \tilde v\cdot\nabla u + \nabla p 
&=& \tau\tilde T\vec\rho ~~~~\mbox{ in } D\nonumber\\
\mbox{div } u &=& 0 ~~~~\mbox{ in } D \label{main1}\\
u = 0 \mbox{ on } \partial D &\mbox{ and }&  
\lim_{x\to\pm\infty}\|u(x,\cdot)\|_{\mathcal{C}^1(\Omega)} = 0.\nonumber
\end{eqnarray}
\begin{theorem}\label{th_uniform_H2}
There exists the unique $u\in H^3\cap \mathcal{C}^{2,\alpha}(D)$ 
solving (\ref{main1}) with some $p\in H^2_{loc}(D)$. Moreover:
\begin{itemize}
\item[(i)] when $d=0$ then $u$ satisfies: 
$\|u\|_{L^\infty(D)}\leq C \|\nabla \tilde T\|_{L^2(R_a)},$
\item[(ii)] when $d=1$ then we have:
$$ \|u\|_{L^\infty(D)}\leq \frac{2 C_P}{\nu\sqrt{\pi\nu}} \left(|\vec\rho|C_{PW}+
 \left(\fint_\Omega |\vec\rho\cdot(0,\tilde x)|^2\right)^{1/2}\right)
\|\tilde v\|_{L^\infty(D)}^{1/2}
\|\nabla \tilde T\|_{L^2(D)} + C\|\nabla \tilde T\|_{L^2(D)}.$$
\end{itemize}
Above, $C_P$ and $C_{PW}$ denote, respectively, the Poincar\'e and the
Poincar\'e-Wirtinger constants of $\Omega$, while the constant $C$ 
is independent of $c, \tau, a, \tilde T$ or $\tilde v$.
\end{theorem}
\begin{proof}
{\em 1. The bound by ${\nabla\tilde T}$.} Define the quantity:
$$L=\left(\fint_\Omega |\vec\rho\cdot(0,\tilde x)|^2\right)^{1/2}$$
and consider the following vector field with boundedly supported gradient:
$$q(x,\tilde x) = \vec\rho\cdot e_1 \int_0^x\fint_\Omega \tilde T(s,\cdot) ~\mbox{d}s
+ \vec\rho\cdot (0,\tilde x) \fint_\Omega \tilde T(x,\cdot).$$
By an easy calculation we see that:
$$ \tilde T(x,\tilde x)\vec\rho - \nabla q(x,\tilde x) =
\left(\tilde T(x,\tilde x) - \fint_\Omega \tilde T(x,\cdot)\right)\vec\rho
- \vec\rho\cdot (0,\tilde x) \fint_\Omega \frac{\partial}{\partial x}\tilde
T(x,\cdot) e_1$$
and therefore:
\begin{equation}\label{quattro}
\|\tilde T\vec\rho - \nabla q\|_{L^2(D)} \leq \left( |\vec\rho| C_{PW} + L\right)
\|\nabla\tilde T\|_{L^2(D)}.
\end{equation}
Recall that the Poincar\'e-Wirtinger constant $C_{PW}$ on $\Omega$ is 
the inverse of the first nonzero eigenvalue of the related Neumann problem.
By a mollification argument, we may also assume that $q\in\mathcal{C}^2(D)$
and that (\ref{quattro}) is still satisfied.

\medskip

{\em 2. Existence of a weak solution.} Following the Galerkin method, 
define:
$$V= cl_{H^1(D)}\left\{u\in\mathcal{C}^\infty_c(D,\mathbb{R}^3), 
~~ \mbox{div } u = 0\right\}.$$
Clearly, $V$ is a Hilbert space with the scalar product
$\langle u,w\rangle_V = \int_{D}\nabla u :\nabla w$.
The norms $\|u\|_V:= \langle u, u\rangle_V^{1/2}$ and $\|u\|_{H^1(D)}$ are
equivalent in $V$, in virtue of the Poincar\'e inequality in $\Omega$,
which yields: $\|u\|_{L^2(D)}\leq C_P \|u\|_V$.

Since $V$ is a subspace of $H^1(D)$, it is also separable and hence it admits a 
Hilbert (orthonormal) basis $\{\psi_n\}_{n=1}^\infty\in \mathcal{C}_c^\infty (D)\cap V$.
For each $n$, let $V_n= \mbox{span }\{\psi_1\ldots \psi_n\}$ and let 
$P_n:V_n\longrightarrow V_n$ be given by:
$$P_n(u) = \sum_{i=1}^n\left\{\nu \langle u,\psi_i\rangle_V 
- c\int_D u_x \psi_i + \tau d \int_D(\tilde v\cdot\nabla u) \psi_i
- \tau\int_D (\tilde T\vec\rho - \nabla q)\psi_i\right\} \psi_i.$$
The operator $P_n$ is continuous and it satisfies:
\begin{equation*}
\begin{split}
\langle P_n(u),&u \rangle_V  = \nu\|u\|_{V}^2 - c\int_D u_x u 
+ \tau d \int_D(\tilde v\cdot \nabla u) u - \tau \int_D (\tilde T\vec\rho -\nabla q)u\\
& \geq \nu \|u\|_V^2 - C_{P}\|\tilde T\vec\rho - \nabla q\|_{L^2(D)} \|u\|_{V} > 0 
\quad \mbox{ when }
\|u\|_V = \frac{2C_P}{\nu} \|\tilde T\vec\rho - \nabla q\|_{L^2(D)},
\end{split}
\end{equation*}
where we used $2\int_D u_x u = \int_D (|u|^2)_x = 0$ and the nullity of
the trilinear term.

By Lemma 2.1.4 in \cite{T}, there exists $u_n\in V_n$, bounded in $V$ by
the above quantity and solving: $P_n(u_n) = 0$. 
Since $V$ is reflexive, it follows that $\{u_n\}$ converges weakly 
(up to a subsequence) to some $u\in V$ such that:
\begin{equation}\label{uno}
\forall w\in V \qquad \nu\int_D \nabla u : \nabla w - c\int_D u_x w 
+ \tau d\int_D (\tilde v\cdot\nabla u) w - \tau\int_D \tilde T \vec\rho w = 0.
\end{equation}
This identity follows first for $w=\psi_n$, and then 
by the density of the linear combinations of $\{\psi_n\}$ in $V$. 
Taking $w=u$ and using (\ref{quattro}) we obtain:
\begin{equation}\label{due}
\|u\|_V\leq \frac{C_P}{\nu}
\left(|\vec\rho| C_{PW} +L \right) \|\nabla\tilde T\|_{L^2(D)}.
\end{equation}

\medskip

{\em 3. Regularity.} 
Recall that by de Rham's theorem 
(see, for example, Proposition 1.1.1 in \cite{T}),
a necessary and sufficient condition for a distribution field $v \in\mathcal{D}'$
that $v=\nabla p$ for some $p\in\mathcal{D}'$ is that $\langle v, w\rangle_V=0$
for all $w\in V$.
Hence, (\ref{uno}) implies the first equality in (\ref{main1})
in the weak sense. By the standard regularity theory \cite{ADN, T} and in view
of (\ref{due}) we may deduce now that the same equality holds in the classical 
sense and that $u\in H^3(D)$, $\nabla p \in H^1(D)$ (since 
$\nabla \tilde T \in L^2(D)$). 
Therefore $u\in\mathcal{C}^{1,\alpha}(D)$ (for $\alpha < 1/4$) and
the boundary conditions in (\ref{main1}) follow, together with the asymptotic
conditions as $|x|\to\infty$, in view of Lemma \ref{lem_pomoc}. 
Next, recalling that $\tilde T \in \mathcal{C}^{1,\alpha}(D)$, 
the potential theory \cite{LU} employed to the localized problem 
and the classical Schauder estimates give that
$u\in\mathcal{C}^{2,\alpha}(D)$. 

\medskip

{\em 4. The bound on ${cu_x}$.}
Since $\nabla q$ has a bounded support and is $\mathcal{C}^1$, 
thus $\nabla (p-q)\in H^1(D)$.  Consequently:
$$\int_D u_x  \nabla (p-q) 
= \int_D (u\nabla (p-q))_x - \int_D \mbox{div } (u (p-q)_x)= 0,$$
because in view of $u\in H^3(D)$ and $\nabla (p-q)\in H^1(D)$ one has:
$$\lim_{|x|\to\infty} \left(\int_\Omega |u\nabla(p-q)|(x,\cdot)
+ \int_\Omega |u^1 (p-q)_x|(x,\cdot)\right) = 0.$$
Integrating the first equality in (\ref{main1}) against $cu_x$ on $D$ we obtain:
\begin{equation*}
\begin{split}
\|cu_x\|_{L^2(D)}^2 & = \nu c\int_D\nabla u : \nabla u_x
- \int_D c u_x (\tau\tilde T\vec\rho - \tau\nabla q) 
+ \tau d \int_D cu_x (\tilde v\cdot\nabla u)\\
& \leq \|cu_x\|_{L^2(D)} \left(\|\tilde T\vec\rho - \nabla q\|_{L^2(D)}
+ d\|\tilde v\cdot\nabla u\|_{L^2(D)}\right),
\end{split}
\end{equation*}
where we have once more used that $u\in H^3(D)$. Therefore, by (\ref{quattro}):
\begin{equation}\label{cinque}
\|cu_x\|_{L^2(D)}\leq 
\left(|\vec\rho|C_{PW}+L\right)
\|\nabla\tilde T\|_{L^2(D)} + d \|\tilde v\cdot\nabla u\|_{L^2(D)}.
\end{equation}

\medskip

{\em 5. The bound and uniqueness for ${d=0}$.} 
It is now easy to conclude the proof,  when
$d=0$. Using  the standard elliptic estimates for the Stokes system (\ref{stokes}),
following from the theory in \cite{ADN}, 
and the Sobolev interpolation inequality, it follows that:
$$\|u\|_{L^\infty(D)}\leq C\|u\|_{H^2(D)}^{1/2}\|u\|_{H^1(D)}^{1/2} 
\leq C \|\nabla\tilde T\|_{L^2(D)},$$
in virtue of (\ref{quattro}), (\ref{due}) and (\ref{cinque}). 
The constant $C$ is uniform and depends only on the geometry of $\Omega$, 
and the constants $|\vec\rho|$ and $\nu$.
Uniqueness of $u$ also follows from the above bound.

\medskip

{\em 6. The case of ${d=1}$.} 
Denote: $g= cu_x - \tau \tilde v\cdot\nabla u + (\tau\tilde T\vec\rho - \tau\nabla q).$
By (\ref{quattro}) and (\ref{cinque}) we obtain that:
$$\|g\|_{L^2(D)}\leq 2\|\tilde v\|_{L^\infty(D)}\|\nabla u\|_{L^2(D)}
+ 2 \left(|\vec\rho| C_{PW} + L\right)
\|\nabla \tilde T\|_{L^2(D)}.$$
Therefore, Theorem \ref{thour_bound} implies:
\begin{equation*}
\begin{split}
\|u\|_{L^\infty(D)} \leq & \frac{2}{\sqrt{\pi\nu}}
\|\tilde v\|_{L^\infty(D)}^{1/2}\|\nabla u\|_{L^2(D)} \\
& + \frac{2}{\sqrt{\pi\nu}}
\left(|\vec\rho| C_{PW} +L\right)^{1/2}
\|\nabla\tilde T\|_{L^2(D)}^{1/2}\|\nabla
u\|_{L^2(D)}^{1/2}  + C_\Omega \|\nabla u\|_{L^2(D)},
\end{split}
\end{equation*}
which by (\ref{due}) establishes the result.
\end{proof}

\section{The uniform bounds and existence of traveling waves}

In this section we prove the uniform bounds on solutions to system
(\ref{trav_wav}), and then establish existence of a traveling wave in
(\ref{nsb}) by a Leray-Schauder degree argument.

Given $c\in\mathbb{R}$, $\tau\in [0,1]$, 
a divergence-free vector field 
$v\in\mathcal{C}^{1,\alpha}(R_a)$ and $Z\in\mathcal{C}^{1,\alpha}(R_a)$
consider first the reaction-advection-diffusion problem:
\begin{eqnarray}
-cT_x  - \Delta T + \tau v\cdot \nabla T & = &\tau f(Z) ~~~~\mbox{ in } R_a\nonumber\\
T(-a, \tilde x) = 1, ~~T(a,\tilde x) = 0 &\mbox{for}& \tilde x\in\Omega \label{rd}\\
\frac{\partial T}{\partial \vec n} (x,\tilde x)= 0 &\mbox{for}& x\in [-a,a] 
\mbox{ and }\tilde x \in\partial\Omega, \nonumber
\end{eqnarray}
together with the following normalization condition, whose eventual role is 
to single out a correct approximation of the traveling wave in $T$, 
in the moving frame which chooses to have $f(T(x,\cdot))=0$ for $x\geq 0$:
\begin{equation}\label{norm}
\max\big\{T(x,\tilde x); ~~ x\in [0,a], ~ \tilde x\in\Omega\big\} = \theta_0.
\end{equation}
We now recall the bounds on solutions to the above problems, 
proved in \cite{BCR} and used in \cite{CLR, Le}. The
right hand sides of (iii), (iv) and (v) follow by re-examining the proofs.
\begin{theorem}\label{th_rd}
Let $T=Z\in\mathcal{C}^{1,\alpha}(R_a)$ satisfy (\ref{rd}) and (\ref{norm}). 
Then one has:
\begin{itemize}
\item[(i)] $\displaystyle{T(x,\tilde x) \in [0,1]}$ for all $(x,\tilde x)\in R_a$,
\item[(ii)] $\displaystyle{T(x,\tilde x) \leq \theta_0}$ for all  
$x>0, \tilde x\in\Omega$, 
\item[(iii)] $\displaystyle{|c|\leq \|v\|_{L^\infty (R_a)} 
+ 2\|f'\|_{L^\infty([0,1])}^{1/2}}$,
\item[(iv)] $\displaystyle{\|\nabla T\|_{L^2(R_a)}^2\leq |\Omega| \left( 
\frac{7}{2} \|c-v^1\|_{L^\infty (R_a)} + \frac{1}{a}\right)}$,
\item[(v)] $\displaystyle{\int_{R_a} f(T) \leq 
 |\Omega| \left(4 \|c-v^1\|_{L^\infty (R_a)} + \frac{1}{a}\right)}$.
\end{itemize}
\end{theorem}

Given $T\in\mathcal{C}^{1,\alpha}(R_a)$ satisfying boundary conditions as in
(\ref{rd}), we will consider its boundedly supported
$\mathcal{C}^{1,\alpha}(D)$ extension:
\begin{equation}\label{Ttilde}
\tilde T(x,\tilde x) = \left\{\begin{array}{ll}
T(x,\tilde x)  & \mbox{ for } x \in [-a,a]\\
 \phi(|x|-a) \cdot \big(2-T(-2a-x,\tilde x)\big) & \mbox{ for } x < -a\\
- \phi(|x|-a) \cdot T(2a - x,\tilde x) & \mbox{ for } x > a.
\end{array}\right.
\end{equation}
Here $\phi\in\mathcal{C}^\infty (\mathbb{R},[0,1])$ satisfies
$\phi(x)=1$ for $x<1/3$, $\phi(x)=0$ for $x>2/3$ and $\|\nabla
\phi\|_{L^\infty}\leq 4$.

Also, for a divergence-free $v\in\mathcal{C}^{1,\alpha}(R_a)$, let 
$\tilde v\in \mathcal{C}^{1,\alpha}(D)$ be its divergence-free extension,
given in Theorem \ref{vtilde}.
\begin{theorem}\label{cor1}
Let $c\in\mathbb{R}$, $\tau\in [0,1]$, $T=Z\in\mathcal{C}^{1,\alpha}(R_a)$,
$v\in \mathcal{C}^{1,\alpha}(R_a)$ and $u\in\mathcal{C}^{1,\alpha}(D)$ 
satisfy (\ref{main1}), (\ref{rd}), (\ref{norm}), 
with $\tilde T$ and $\tilde v$ defined as above.
Moreover, let $u_{\mid R_a}= v$ and
assume that either $d=0$, or $d=1$ and (\ref{small}) holds.
Then, for large $a$:
\begin{equation*}
|c| + \|T\|_{\mathcal{C}^{1,\alpha}(R_a)} +  \|u\|_{\mathcal{C}^{2,\alpha}(R_a)} + 
\|\nabla T\|_{L^2(R_a)} + \|u\|_{H^3(R_a)} + \int_{R_a} f(T) \leq C,
\end{equation*}
where $C$ is a numeric constant independent on $a$, $\tau$ and 
the estimated quantities. 
\end{theorem}
\begin{proof}
By Theorem \ref{th_rd} (iii) and (iv) we obtain:
$$\|\nabla T\|_{L^2(R_a)}^2 \leq 7|\Omega|\big(\|u\|_{L^\infty(D)}
+ \|f'\|_{L^\infty([0,1])}^{1/2}\big) + \frac{|\Omega|}{a}.$$
On the other hand, the boundary conditions for $T$ imply that, for large $a$:
$$\|\nabla\tilde T\|_{L^2(D)}\leq \sqrt{2} \|\nabla T\|_{L^2(R_a)} + 8
\leq \sqrt{14} |\Omega|^{1/2}\big(\|u\|_{L^\infty(D)}^{1/2}+
\|f'\|_{L^\infty([0,1])}^{1/4}\big) + 9.$$
Consequently, for $d=0$ the uniform bound on $\|u\|_{L^\infty(D)}$ follows 
by Theorem \ref{th_uniform_H2} (i).  

When $d=1$ then Theorem
\ref{th_uniform_H2} (ii) yields the same bound, under condition
(\ref{small}). Indeed, let $\epsilon>0$ be such that the quantity in the
left hand side of (\ref{small}) is strictly smaller than 
$1/\sqrt{1+\epsilon}$. Then:
\begin{equation}\label{inq-main}
\begin{split}
\|u\|_{L^\infty(D)} &\leq \sqrt{1+\epsilon} \frac{2 C_P}{\nu\sqrt{\pi\nu}} 
\left(|\vec\rho|C_{PW}+
 \left(\fint_\Omega |\vec\rho\cdot(0,\tilde x)|^2\right)^{1/2}\right)\times \\
& \qquad\qquad
\times \left(\sqrt{14} |\Omega|^{1/2} \big(\|u\|_{L^\infty(D)}^{1/2} 
+ \|f'\|_{L^\infty([0,1])}^{1/4}\big) + C\right) \|u\|_{L^\infty(D)}^{1/2} \\
& \quad 
+ C\big(\|u\|_{L^\infty(D)}^{1/2} + \|f'\|_{L^\infty([0,1])}^{1/4} + 1\big)\\
&\leq q \|u\|_{L^\infty(D)} + C\big(\|u\|_{L^\infty(D)}^{1/2} +1\big)
\big(\|f'\|_{L^\infty([0,1])}^{1/4} + 1\big),
\end{split}
\end{equation}
for some $q\in(0,1)$. Hence
$\|u\|_{L^\infty(D)}\leq C \big(\|f'\|_{L^\infty([0,1])}^{1/2} + 1\big)$
and so, recalling (\ref{due}) and Theorem \ref{th_uniform_H2}, we have:
\begin{equation}\label{piotr}
|c| + \|u\|_{L^\infty(D)} + \|u\|_{H^1(D)}^2 + \int_{R_a} f(T) 
\leq C \big(\|f'\|_{L^\infty([0,1])}^{1/2} + 1\big).
\end{equation}
In (\ref{inq-main}) and (\ref{piotr}) the constant $C$ is independent of
$a$, $\tau$, the nonlinearity $f$ and the estimated quantities.

The uniform bounds on $\|u\|_{H^2(D)}$, $|c|$, 
$\|\nabla T\|_{L^2(R_a)}$ and $\int_{R_a} f(T)$ follow by
Theorem \ref{th_rd}, (\ref{due}) and (\ref{lapl}). 
Now, the standard local elliptic estimates for the Stokes system 
(\ref{stokes}) (see \cite{ADN}, also compare \cite{Le}) 
imply that:
$$\|u\|_{H^3(D)}\leq C (\|\nabla g\|_{L^2(D)}  + \|u\|_{H^1(D)}).$$
Taking $g= c u_x - \tau\tilde v\cdot \nabla u + \tau\tilde T\vec\rho$, 
we obtain the uniform bound on $\|u\|_{H^3(D)}$.
The bounds on $\|u\|_{\mathcal{C}^{2,\alpha}(D)}$ and 
$\|T\|_{\mathcal{C}^{1,\alpha}(R_a)}$ follow by  H\"older's estimates 
for system (\ref{main1}) and (\ref{rd}) \cite{ADN, GT}.
The proof is done.

We remark that our result does not imply smallness of $C$ in the 
uniform bound. In particular, $C$ depends on the constant $C_{\epsilon,\Omega}$
from Theorem \cite{LLP}, which can be arbitrarily large.
\end{proof}

\bigskip

\noindent We finally have:

{\bf Proof of Theorem \ref{th2}.}
For every sufficiently large $a>0$, consider an operator:
\begin{equation*}
K_a:\mathbb{R}\times\mathcal{C}^{1,\alpha}(R_a)
\times\mathcal{C}_d^{1,\alpha}(R_a) \times [0,1] \longrightarrow
\mathbb{R}\times \mathcal{C}^{1,\alpha}(R_a)\times\mathcal{C}_d^{1,\alpha}(R_a),
\end{equation*}
where $\mathcal{C}_d^{1,\alpha}(R_a)$ stands for the Banach space of the
divergence-free, $\mathcal{C}^{1,\alpha}$ regular vector fields 
on the compact domain $R_a$. Define:
\begin{equation*} K_a(c,Z,v,\tau):= \big(c-\theta_0 + \max\{T(x,\tilde x);~~ x\in [0,a], 
~ \tilde x\in\Omega\}, T, u_{\mid R_a}\big),\nonumber
\end{equation*}
where $T$ is the solution to (\ref{rd}), and $u$ solves (\ref{main1}) are
known also
with $\tilde T$ and $\tilde v$ defined as in (\ref{Ttilde})
and Theorem \ref{vtilde}.

The operator $K_a$ is continuous, compact \cite{GT}
and all its fixed points $(c,T, v)$ such that
$K_a(c,T,v,\tau) = (c,T, v)$ for some $\tau\in [0,1]$ are uniformly bounded, in view
of Theorem \ref{cor1}.
We may now employ the Leray-Schauder degree theory, as in \cite{BCR, CLR, Le},
to obtain the existence of a fixed point of $K_a(\cdot,\cdot,\cdot, 1)$,
since the degree of the map $K_a(\cdot,\cdot,\cdot,0)$ is nonzero.
This fixed point $(c^a, T^a, v^a)$ again satisfies the bounds in Theorem
\ref{cor1}.

By a bootstrap argument we moreover obtain the uniform bound on
$\|T^a\|_{\mathcal{C}^{2,\alpha}(R_{a-1})} $.
One may thus choose a sequence $a_n\to\infty$ such that $c_n:=c^{a_n}$ converges
to some $c\in\mathbb{R}$, and $T_n:=T^{a_n}$, $v_n:=v^{a_n}$ converge in
$\mathcal{C}_{loc}^{2,\alpha}(D)$ to some $T,u\in\mathcal{C}^{2,\alpha}(D)$.
Further, $u \in H^2(D)\cap \mathcal{C}^{2,\alpha}(D)$ and hence 
the first convergence in (\ref{tre.due}) follows. 
Since $\nabla T\in L^2\cap \mathcal{C}^{0,\alpha}(D)$, we obtain
the other convergence in view of Lemma \ref{lem_pomoc}. 

The positivity of the propagation speed $c$ and
the existence of the right and left limits of $T$, together with 
the statement in (ii) follow exactly as in \cite{Le}.
\endproof

\section{A sufficient condition for $\theta_-=1$}

The following lemma improves on the result in \cite{Le}, where a sufficient
condition for the left limit $\theta_-$ of $T$ 
to be $1$ required a cubic bound: $f(T) \leq k [(T-\theta_0)_+]^3$.

\begin{lemma}\label{ththeta-}
In the setting of Theorem \ref{th2}, if moreover the nonlinearity 
satisfies:
\begin{equation}\label{cond_small}
f(T) \leq k [(T-\theta_0)_+]^2 \quad \mbox{ and } \quad
k\left(1 + \|f\|_{L^\infty([0,1])} + \|f'\|_{L^\infty([0,1])}\right)\leq C_\Omega,
\end{equation}
then $\theta_-=1$.
Here $C_\Omega>0$ is a constant, depending only on $\nu$, $\rho$ 
and $\Omega$.
\end{lemma}
\begin{proof}
{\em 1.} 
In the course of the proof, $C$ will denote any positive constant depending 
only on $\nu$, $\rho$  and $\Omega$. 
Integrating the temperature equation in (\ref{trav_wav}) against $T$ and
$\Delta $ on $D$ yields, respectively:
\begin{eqnarray}
{\displaystyle{
\|\nabla T\|^2_{L^2(D)}}}&=& {\displaystyle{\int_D f(T)T - \frac{1}{2}c\theta_-^2 |\Omega|
\leq \int_D f(T)}}\label{p0}\\
\|\Delta T\|_{L^2(D)} &\leq& 
{\displaystyle{\|f(T)\|_{L^2(D)} + \|u\cdot \nabla T\|_{L^2(D)},}}
\label{p1}
\end{eqnarray}
where in (\ref{p1}) we used that $\int_D T_x\Delta T=0$.
The interpolation, H\"older and Sobolev inequalities imply that:
\begin{equation}\label{p2}
\begin{split}
\|u\cdot \nabla T\|_{L^2(D)} &\leq \|u\|_{L^6(D)} \|\nabla T\|_{L^3(D)}
\leq \|u\|_{L^6(D)} \|\nabla T\|_{L^2(D)}^{1/2}\|\nabla T\|_{L^6(D)}^{1/2}\\
&\leq \frac{C}{\epsilon} \|u\|_{H^1(D)}^2 \|\nabla T\|_{L^2(D)}
+ \epsilon \|\nabla T\|_{H^1(D)}.
\end{split}
\end{equation}
Now, taking $\epsilon$ above sufficiently small and introducing (\ref{p1})
and (\ref{p2}) into:
$$ \|\nabla T\|_{H^1(D)}\leq C\left(\|\nabla T\|_{L^2(D)} + \|\Delta T\|_{L^2(D)} \right)$$
we obtain, in view of (\ref{p0}):
\begin{equation*}
\begin{split}
\|\nabla T\|_{H^1(D)} &\leq C\left(\|\nabla T\|_{L^2(D)} + \|f(T)\|_{L^2(D)}
+ \|u\|^2_{H^1(D)}\|\nabla T\|_{L^2(D)}\right)\\
&\leq C\left(1 + \|f\|_{L^\infty([0,1])}^{1/2} + \|u\|^2_{H^1(D)}\right)
\left(\int_Df(T)\right)^{1/2}.
\end{split}
\end{equation*}
By (\ref{piotr}) and convergences established in the proof of Theorem
\ref{th2},
this implies:
\begin{equation}\label{p3}
\|\nabla T\|_{H^1(D)} \leq C\left(1 + \|f\|_{L^\infty}^{1/2} + \|f'\|^{1/2}_{L^\infty}\right)
\left(\int_Df(T)\right)^{1/2}.
\end{equation}

{\em 2.} 
Now, for every $x\in\mathbb{R}$ denote $M(x) = \max_{\tilde x\in\Omega} T(x,\tilde x)$,
$m(x) = \min_{\tilde x\in\Omega} T(x,\tilde x)$ and
notice that $m(x)$ is non increasing. This can be proved for
each $T_n$ on $R_n$, using the maximum principle. 
Passing with $n$ to $\infty$, one obtains the same result in the limit.

We now argue by contradiction. If $\theta_-\leq \theta_0$ then 
$m(x)\leq \theta_0$ for every $x\in\mathbb{R}$ and:
\begin{equation*}
\begin{split}
\int_D [(T-\theta_0)_+]^2 &\leq |\Omega|  \int_{-\infty}^{+\infty} |M(x) -
m(x)|^2 ~\mbox{d}x\\
&\leq 2 |\Omega|  \int_{-\infty}^{+\infty}\|T(x,\cdot) 
- \fint_\Omega T(x,\cdot)\|_{L^\infty(\Omega)}~\mbox{d}x
\leq C \|\nabla T\|^2_{H^1(D)}.
\end{split}
\end{equation*}
Together with (\ref{p3}) and the assumption in (\ref{cond_small}) the above yields:
$$\int_D [(T-\theta_0)_+]^2\leq 
C\left(1 + \|f\|_{L^\infty} + \|f'\|_{L^\infty}\right) k \int_D [(T-\theta_0)_+]^2,$$
which by the assumption on $k$ implies that both sides above must be
zero. Consequently, be $f(T)\equiv 0$ and one can deduce (see \cite{BCR,Le})
that $T\equiv 0$ as well, contradicting the results of Theorem \ref{th2}.
\end{proof}

The condition (\ref{cond_small}) seems to be artificial and we believe that
it can be further relaxed or even omitted altogether, for the wave $(T,u)$ 
obtained in the limiting procedure of Theorem \ref{th2}.


\begin{thebibliography}{999999}
\bibitem{ADN} S. Agmon, A. Douglis and L. Nirenberg, 
\textit{Estimates near the boundary 
for solutions of elliptic partial differential equations satisfying general 
boundary conditions. II}, Comm. Pure Appl. Math.,  {\bf 17}  (1964), 35--92.

\bibitem{B} H. Berestycki, \textit{The influence of advection on
the propagation of fronts in reaction-diffusion equations}, 
Nonlinear PDEs in Condensed Matter and Reactive Flows, NATO Science
Series C, {\bf 569}, H. Berestycki and Y. Pomeau eds, Kluwer, Doordrecht,
2003.

\bibitem{BCR} H. Berestycki, P. Constantin and L. Ryzhik,
  \textit{Non-planar fronts in Boussinesq reactive flows},
 Ann. Inst. H. Poincar\'e Anal. Non Lin\'eaire,  {\bf 23}  (2006),  
no. 4, 407--437. 

\bibitem{CKR} P. Constantin, A. Kiselev and L. Ryzhik, \textit{
  Fronts in reactive convection: bounds, stability and instability}
  Comm. Pure Appl. Math., {\bf 56} (2003), 1781--1803.

\bibitem{CLR} P. Constantin, M. Lewicka and L. Ryzhik, 
\textit{Traveling waves in 2D reactive Boussinesq systems 
with no-slip boundary conditions}, Nonlinearity, {\bf 19} (2006), 2605--2615.

\bibitem{G} G.P. Galdi, \textit{An introduction to the mathematical theory 
of the Navier-Stokes equations}, Springer Tracts in Natural Philosophy, {\bf 38. 39},
Springer, 1998.

\bibitem{GT} D. Gilbarg and N. Trudinger, \textit{Elliptic partial
    differential equations of second order}, Springer-Verlag, Berlin, 2001.

\bibitem{LU} O.A. Ladyzhenskaya and N.N. Uraltseva, \textit{Linear and
    quasilinear elliptic equations}, Translated from the Russian by Scripta
  Technica, Inc. Translation editor: Leon Ehrenpreis Academic Press, 1968.

\bibitem{La} O.A. Ladyzhenskaya, \textit{On unique solvability 
of three-dimensional Cauchy problem for the Navier--Stokes equations under 
the axial symmetry}, Za. Nauchn. Sem. LOMI 7 (1968) 
 155--177 (in Russian).

\bibitem{Le} M. Lewicka, \textit{Existence of traveling waves in the 
Stokes-Boussinesq system for reactive flows}, J. Diff. Equations, {\bf 237} 
(2007), no. 2, 343--371.

\bibitem{LLP} J-G. Liu, J. Liu and R.L. Pego, \textit{Stability and
    convergence of efficient Navier-Stokes solvers via a commutator estimate},  
Comm. Pure Appl. Math.  {\bf 60}  (2007),  no. 10, 1443--1487.

\bibitem{Mu1} P.B. Mucha, \textit{ On a problem for the Navier-Stokes
    equations with the infinite
 Dirichlet integral},  Z. Angew. Math. Phys.  {\bf 56}  (2005),  no. 3, 439--452.

\bibitem{Mu2} P.B. Mucha, \textit{On cylindrical symmetric flows through 
pipe-like domains}, J. Differential Equations  {\bf 201}  (2004),  no. 2, 304--323.

\bibitem{R} J.M. Roquejoffre,  \textit{Eventual monotonicity and convergence to traveling 
   fronts for the solutions of parabolic equations in cylinders},
   Ann. Inst. Henri Poincare, {\bf 14} (1997), no 4, 499--552.

\bibitem{S} H. Sohr, \textit{The Navier-Stokes equations. An elementary
    functional analytic approach}, Birkh\"auser Verlag, Basel, 2001.

\bibitem{T} R. Temam, \textit{Navier-Stokes equations. Theory and numerical
    analysis.} AMS Chelsea Publishing, Providence, RI, 2001.

\bibitem{TV} R. Texier-Picard and V. Volpert, 
   \textit{Problemes de reaction-diffusion-convection dans des cylindres non bornes},
   C. R. Acad. Sci. Paris Sr. I Math., {\bf 333} (2001), 1077--1082.

\bibitem{Yu} M.R. Ukhovskij and V.I. Yudovich, 
\textit{Axially symmetric motions of prefect and 
viscous fluids filling all space}, Prihl. Math. Mekh. {\bf 32} (1968), 59--69 (in Russian).

\bibitem{VR1} N. Vladimirova  and R. Rosner, 
   \textit{Model flames in the Boussinesq limit: the efects of feedback},
   Phys. Rev. E., {\bf 67} (2003), 066305.                       

\bibitem{VR2}  N. Vladimirova  and R. Rosner, 
   \textit{Model flames in the Boussinesq limit: the case of pulsating fronts},
Phys. Rev. E., {\bf 71} (2005), 067303.   

\bibitem{Xie} W. Xie, \textit{A sharp pointwise bound for functions with 
$L\sp 2$-Laplacians and zero boundary values of arbitrary three-dimensional
domains}, Indiana Univ. Math. J.  {\bf 40}  (1991),  no. 4, 1185--1192.

\bibitem{XieStokes}  W. Xie, \textit{On a three-norm inequality for the Stokes 
operator in nonsmooth domains},  
The Navier-Stokes equations II---theory and numerical methods, 
Springer Lecture Notes in Math., {\bf 1530} (1992),  310--315.

\bibitem{Xin} J. Xin, \textit{Front propagation in heterogeneous media},
   SIAM Review {\bf 42} (2000), no 2, 161--230.

\bibitem{Z} Ya.B. Zeldovich, G.I. Barenblatt, V.B. Librovich and G.M.
Makhviladze, \textit{The Mathematical Theory of Combustion and
Explosions}, Consultants Bureau, New York, 1985.
\end{thebibliography}
\end{document}